\documentclass[12pt, reqno]{amsart}
\usepackage{amssymb}
\usepackage{amsfonts}
\usepackage{amssymb,latexsym}
\usepackage{enumerate}
\usepackage{url}
\usepackage{mathrsfs}
\usepackage{bbm}

\allowdisplaybreaks

\makeatletter
\@namedef{subjclassname@2010}{%
  \textup{2010} Mathematics Subject Classification}
\makeatother

  [2002/01/19 v2.2g %
    AMS font definitions%
  ]
\DeclareFontFamily{U}{euf}{}
\DeclareFontShape{U}{euf}{m}{n}{%
  <5><6><7><8><9>gen*eufm%
  <10><10.95><12><14.4><17.28><20.74><24.88>eufm10%
  }{}
\DeclareFontShape{U}{euf}{b}{n}{%
  <5><6><7><8><9>gen*eufb%
  <10><10.95><12><14.4><17.28><20.74><24.88>eufb10%
  }{}

  [2002/01/19 v2.2g %
    AMS font definitions%
  ]
\DeclareFontFamily{U}{msb}{}
\DeclareFontShape{U}{msb}{m}{n}{%
  <5><6><7><8><9>gen*msbm%
  <10><10.95><12><14.4><17.28><20.74><24.88>msbm10%
  }{}

  [2002/01/19 v2.2g %
    AMS font definitions%
  ]
\DeclareFontFamily{U}{msa}{}
\DeclareFontShape{U}{msa}{m}{n}{%
  <5><6><7><8><9>gen*msam%
  <10><10.95><12><14.4><17.28><20.74><24.88>msam10%
  }{}

\newtheorem{theorem}{Theorem}[section]

\newtheorem{corollary}[theorem]{Corollary}

\theoremstyle{definition}
\newtheorem{definition}[theorem]{Definition}

\newtheorem{remark}[theorem]{Remark}

\numberwithin{equation}{section} 

\begin{document}

\title[Ramanujan's identity for Bernoulli-Carlitz numbers]{An analogue of Ramanujan's identity for Bernoulli-Carlitz numbers}

\author{Su Hu}
\address{Department of Mathematics, South China University of Technology, Guangzhou, Guangdong 510640, China}
\email{mahusu@scut.edu.cn}

\author{Min-Soo Kim}
\address{Department of Mathematics Education, Kyungnam University, Changwon, Gyeongnam 51767, Republic of Korea}
\email{mskim@kyungnam.ac.kr}


\begin{abstract}
In his second notebook, Ramanujan discovered the following identity for the special values of
$\zeta(s)$ at the odd positive integers
		\begin{equation*}\begin{aligned}
				\alpha^{-m}\,&\left\{\dfrac{1}{2}\,\zeta(2m + 1) + \sum_{n = 1}^{\infty}\dfrac{n^{-2m - 1}}{e^{2\alpha n} - 1}\right\} \\
				&\quad\quad\quad\quad\quad\quad\quad-(- \beta)^{-m}\,\left\{\dfrac{1}{2}\,\zeta(2m + 1) + \sum_{n = 1}^{\infty}\dfrac{n^{-2m - 1}}{e^{2\beta n} - 1}\right\}\nonumber\\
	&=2^{2m}\sum_{k = 0}^{m + 1}\dfrac{\left(-1\right)^{k-1}B_{2k}\,B_{2m - 2k+2}}{\left(2k\right)!\left(2m -2k+2\right)!}\,\alpha^{m - k + 1}\beta^k,\end{aligned}
	\end{equation*}
where $ \alpha $ and $ \beta $ are positive numbers such that $ \alpha\beta = \pi^2 $ and $ m $ is a positive integer.
As shown by Berndt in the viewpoint of general transformation
 of analytic Eisenstein series, it is a natural companion of  Euler's famous formula for even zeta values.

In this note, we prove an analogue of the above Ramanujan's identity in the functions fields setting, which involves the Bernoulli-Carlitz numbers.
\end{abstract}

\maketitle

\section{Introduction}
Let 
$$
\zeta(s)=\sum_{n=1}^{\infty}\frac{1}{n^{s}}, ~~\textrm{Re}(s)>1
$$
be the Riemann zeta function. It can be analytically continued to the whole complex plane
except for a single pole at $s=1$.

Around 1742, Euler obtained the special values of $\zeta(s)$ at the positive even integers
\begin{equation}\label{(1.1)}
\zeta(2m)=(2\pi)^{2m}\frac{(-1)^{m}B_{2m}}{2(2m)!},
\end{equation}
where $B_{2m}$ are the Bernoulli numbers defined by the generating function
\begin{equation}\label{Bernoulli} 
\frac{z}{e^{z}-1}=\sum_{m=0}^{\infty}B_{m}\frac{z^{m}}{m!}.
\end{equation}
Then Ramanujan discovered the following formula for the special values of
$\zeta(s)$ at the odd positive integers (see Entry 21(i) in Ramanujan's second notebook 
\cite[p. 173]{Ramanujan 2} or his lost notebook \cite[p. 319--320, formula (28)]{Ramanujan 3}).

\begin{theorem}[Ramanujan's formula for $ \zeta(2m+1) $]\label{Ramanujan}
	If $ \alpha $ and $ \beta $ are positive numbers such that $ \alpha\beta = \pi^2 $ and if $ m $ is a positive integer, then we have
	\begin{equation}\label{(1.2)}
	\begin{aligned}
		\alpha^{-m}\,&\left\{\dfrac{1}{2}\,\zeta(2m + 1) + \sum_{n = 1}^{\infty}\dfrac{n^{-2m - 1}}{e^{2\alpha n} - 1}\right\} \\
		&\quad\quad\quad\quad\quad\quad\quad-(- \beta)^{-m}\,\left\{\dfrac{1}{2}\,\zeta(2m + 1) + \sum_{n = 1}^{\infty}\dfrac{n^{-2m - 1}}{e^{2\beta n} - 1}\right\} \\
	&=2^{2m}\sum_{k = 0}^{m + 1}\dfrac{\left(-1\right)^{k-1}B_{2k}\,B_{2m - 2k+2}}{\left(2k\right)!\left(2m -2k+2\right)!}\,\alpha^{m - k + 1}\beta^k. 
	\end{aligned}
   \end{equation}	
\end{theorem}

In 1925, Marulkar \cite{SL} provided the first proof of equation (\ref{(1.2)}), though he was not aware that this formula had already appeared in Ramanujan’s notebooks. Later, (\ref{(1.2)}) was rediscovered and studied by Grosswald in \cite{Gross1, Gross2}.
In 1977, Berndt \cite[Theorem 2.2]{Berndt} demonstrated that both Euler's formula (\ref{(1.1)}) and Ramanujan's formula (\ref{(1.2)}) are special cases of a general transformation formula for analytic Eisenstein series, establishing them as natural companions. In 2011, Gun, Murty, and Rath \cite{GMR} offered a modern interpretation of (\ref{(1.2)}), showing that it encodes fundamental transformation properties of Eisenstein series on the full modular group and their Eichler integrals. This approach was further extended by Berndt and Straub \cite[Section 5]{BS} to weight $2k+1$ Eisenstein series of level 2 via secant Dirichlet series.
In 2024, as a complement to the work of Berndt and Straub \cite{BS}, Dixit \cite{Dixit 2024} provided a survey of recent developments and generalizations of Ramanujan's identity (\ref{(1.2)}), covering extensions in settings such as Lambert series, Koshliakov zeta functions, higher Herglotz functions, and non-holomorphic Eisenstein series. The survey also discussed connections to modular forms, number theory, and other mathematical areas that have emerged in the last decade.

Recently, Chavan \cite{Chavan} presented an elementary proof of (\ref{(1.2)}), based on a Mittag-Leffler type expansion of the hyperbolic cotangent function:
\begin{equation}\label{ML2}
\coth(\pi x) = \frac{1}{\pi x} + \frac{2x}{\pi} \sum_{k=1}^{\infty} \frac{1}{x^2 + k^2},
\end{equation}
for $x \in \mathbb{C} \setminus \{0\}$, and Euler's identity for even zeta values (\ref{(1.1)}).

In \cite{DKK}, the authors generalized Ramanujan's identity (\ref{(1.2)}) and situated it within a broader analytic framework. Specifically,
Ramanujan's identity is derived as a special case of their master identity (Theorem 2.5) when $a = -2m - 1$,
 this master identity provides a unified transformation for the Lambert series $\sum_{n=1}^\infty \sigma_a(n)e^{-ny}$, covering both modular and non-modular cases,
and  for an odd integer $a$, the master identity reduces to classical modular transformations (e.g., Eisenstein series on $\mathrm{SL}_2(\mathbb{Z})$), which are connected to Ramanujan's formula.
Furthermore, the authors also derived a new companion formula (Corollary 2.13) involving even zeta values $\zeta(2m)$, contrasting with Ramanujan's result for odd zeta values.
Thus, the paper not only recovers Ramanujan's result but also extends it to a broader class of series and provides new identities in the same spirit.

In this note, we shall investigate the above Ramanujan's identity in the functions fields setting. In this situation, an analogy of the Riemann zeta function $\zeta(s)$ was introduced by Goss \cite{Goss} in 1979.
Let $q=p^{k}$ be a power of a prime number $p$, $\mathbb{F}_{q}$ be the finite field with $q$ elements and $\mathbb{F}_{q}^{*}$ be its multiplicative group.
Let $A=\mathbb{F}_{q}[T]$ be the polynomial ring with one variable
over the finite field $\mathbb{F}_{q}$ and $K=\mathbb{F}_{q}(T)$ be the rational function field. Let $K_{\infty}=\mathbb{F}_{q}((\frac{1}{T}))$ be the completion of $K$ at the infinite place $\infty=\left(\frac{1}{T}\right)$,
$K_{\infty}^{*}=K_{\infty}\setminus\{0\}$ be the multiplicative group of the field $K_{\infty}$ and $\mathbb{C}_{\infty}$ be the completion of an algebraic closure of $K_{\infty}$.
Let $A_{+}$ be the set of monic polynomials in $A,$ which is an analogue of the set of positive integers.
Goss considered the infinite series
\begin{equation}\label{(1.3)}
\zeta_{\infty}(s)=\sum_{a\in A_{+}}a^{-s}
\end{equation}
for $s\in\mathbb{S}=K_{\infty}^{*}\times\mathbb{Z}_{p},$ where $\mathbb{Z}_{p}$ is the ring of $p$-adic integers,
and he also showed that it is an analytic function on $\mathbb{S}$. For details, we refer to Goss' original
article \cite{Goss}.

Let 
$$[j] =T^{q^j}-T,\quad D_{j}=[j]{[j-1]}^q \cdots {[1]}^{q^{j-1}}$$
with $D_{0}=1$. 
Let \( C \) be the Carlitz module  and the Carlitz exponential $e_{C}(z)$ associated with $C$ is defined by
\begin{equation}\label{1.6+}
e_{C}(z)=\sum_{j=0}^\infty\frac{z^{q^j}}{D_j}.
\end{equation}
It is an entire $\mathbb{F}_{q}$-linear function and has the following infinite product expansion
\begin{equation}\label{(1.6)}
e_{C}(z)=e_{\tilde{\pi}A}(z)=z\prod_{0\neq \lambda\in\tilde{\pi}A}\left(1-\frac{z}{\lambda}\right),
\end{equation}
where $\tilde{\pi}\in\mathbb{C}_{\infty}$ is an analogy of  the classical $2\pi i~(i=\sqrt{-1})$, as it is 
a fundamental period of Carlitz's exponential $e_{C}(z)$.
In 1941, Wade~\cite{Wade} proved that $\tilde{\pi}$ is transcendental over $K$.

For a nonnegative integer $m$ with $q$-ary expansion
$m=\sum_{j=0}^{w}c_{j}q^{j}~(0\leq c_{j}<q),$
the Carlitz factorial $\Gamma_m$ is 
\begin{equation}\label{factorial}
\Gamma_m=\prod_{j=0}^{w}D_{j}^{c_j}.
\end{equation}
Then in analogy with the definition of the classical Bernoulli numbers (\ref{Bernoulli}), 
the Bernoulli-Carlitz numbers $BC_{m}\in K~ (n=0,1,2,\ldots)$
are defined by the generating function
\begin{equation}\label{(1.5)}
\frac{z}{e_{C}(z)}=\sum_{m=0}^{\infty}\frac{BC_{m}}{\Gamma_{m}}{z^{m}}.
\end{equation}
As in the classical case, the Bernoulli-Carlitz numbers  have many deep applications in the arithmetic of function fields, such as Kummer's and  Herbrand's criterions for the class
groups of cyclotomic function fields (see \cite[Section 9.2]{Go3}, \cite[Section 5.2]{Thakurbook} or \cite{Gekeler1, Gekeler2}).
In 1930s, Carlitz \cite{Car1, Car2} showed an analogue of Euler's formula (\ref{(1.1)}), that is, the special values of  $\zeta_{\infty}(s)$
at the “even" positive integers $m=(q-1)j~(j\in\mathbb{N})$ are given by 
\begin{equation}\label{(1.4)}
\zeta_{\infty}(m)=\tilde{\pi}^{m}\frac{BC_{m}}{\Gamma_{m}}.
\end{equation}
For details on the special values of $\zeta_{\infty}(s)$, we refer to Thakur's work \cite{Thakur}.

Now suppose $q$ is a power of an odd  prime and by understanding that $BC_{0}=\Gamma_{0}=1$. In this note, we shall prove the following function fields analogue of Ramanujan's identity (\ref{(1.2)}).
\begin{theorem}\label{main}
Let $e_{C}(z)$ be Carlitz's exponential function (\ref{1.6+}), $\Gamma_{m}$ be Carlitz's  factorial  (\ref{factorial}) and $BC_{m}$ be the Bernoulli-Carlitz numbers (\ref{(1.5)}). If $ \alpha $ and $ \beta $ are elements in $\mathbb{C}_{\infty}$ such that $ \alpha\beta = \tilde{\pi}^2 $ and if $ j $ is a positive integer, then we have
\begin{align}
\alpha^{-(\frac{q-1}{2})(j-1)-(q-2)}&\sum_{a\in A_{+}}\frac{a^{-(q-1)j-(q-2)}}{e_{C}(\alpha a)}+\beta^{-(\frac{q-1}{2})(j-1)-(q-2)}\sum_{a\in A_{+}}\frac{a^{-(q-1)j-(q-2)}}{e_{C}(\beta a)} \nonumber\\
&=\sum_{k=0}^{j+1}\frac{BC_{(q-1)k}}{\Gamma_{(q-1)k}}\frac{BC_{(q-1)(j-k+1)}}{\Gamma_{(q-1)(j-k+1)}}\alpha^{(\frac{q-1}{2})k}\beta^{(\frac{q-1}{2})(j-k+1)}. \label{1.7}\end{align}
\end{theorem}
\begin{remark}
By comparing (\ref{1.7}) with (\ref{(1.2)}), we see that  the terms 
$$\sum_{a\in A_{+}}\frac{a^{-(q-1)j-(q-2)}}{e_{C}(\alpha a)}\quad\text{and}\quad\sum_{a\in A_{+}}\frac{a^{-(q-1)j-(q-2)}}{e_{C}(\beta a)}$$
are analogies of
$$\sum_{n = 1}^{\infty}\dfrac{n^{-2m - 1}}{e^{2\alpha n} - 1} \quad\text{and}\quad \sum_{n = 1}^{\infty}\dfrac{n^{-2m - 1}}{e^{2\beta n} - 1},$$
respectively.
\end{remark}

The subsequent sections of this paper are organized as follows. In Section 2, we provide some remarks that connect the above result with the theory of modular forms over function fields, including relationships between Petrov's $A$-expansions and Hamahata's construction of Lambert series. In Section 3, we present the proof of Theorem \ref{main}.

\section{On its relations between modular forms in function fields}
In this section, as suggested by the referee, we give some remarks which relate Theorem \ref{main} with the theory of modular forms in the function fields.

\subsection{Petrov's $A$-expansions} In the classical setting, let the upper half-plane be $\mathbb{H} = \{ \tau \in \mathbb{C} \mid \operatorname{Im}(\tau) > 0 \}$. The modular group $\mathrm{SL}_2(\mathbb{Z})$ acts on $\mathbb{H}$ by fractional linear transformations: for a matrix $\begin{pmatrix} a & b \\ c & d \end{pmatrix} \in \mathrm{SL}_2(\mathbb{Z})$, the action is defined by $\tau \mapsto \frac{a\tau + b}{c\tau + d}$. 
For any odd integer $n \geq 1$, recall that the Lambert series is defined as
\begin{equation}\label{Gn}
G_n(x) = \sum_{m=1}^{\infty} \frac{m^{-n} x^m}{1 - x^m}, \quad |x| < 1.
\end{equation}
By using $G_{n}(x)$, Ramanujan's identity (\ref{(1.2)}) can be wrtten as 
\begin{equation}\label{(1.2+)}
	\begin{aligned}
		\alpha^{-m}\,&\left\{\dfrac{1}{2}\,\zeta(2m + 1) + G_{2m+1}\left(e^{-2\alpha}\right)\right\} \\
		&\quad\quad\quad\quad\quad\quad\quad-(- \beta)^{-m}\,\left\{\dfrac{1}{2}\,\zeta(2m + 1) +  G_{2m+1}\left(e^{-2\beta}\right)\right\} \\
	&=2^{2m}\sum_{k = 0}^{m + 1}\dfrac{\left(-1\right)^{k-1}B_{2k}\,B_{2m - 2k+2}}{\left(2k\right)!\left(2m -2k+2\right)!}\,\alpha^{m - k + 1}\beta^k. 
	\end{aligned}
   \end{equation}

Let
\[
\mathcal{G}_{2k}(z) := \sum_{(0,0)\neq(m,n)} \frac{1}{(mz + n)^{2k}}
\]  
and 
\[
\mathcal{E}_{2k}(z) := \mathcal{G}_{2k}(z)/2\zeta(2k)
\]  
be the Eisenstein series and its normalized version, respectively. From the Fourier expansion of $\mathcal{G}_{2k}(z)$,
we have the following Lambert expansion for \(\mathcal{E}_{2k}(z)\):
\[
\mathcal{E}_{2k}(z) = 1 + \frac{2}{\zeta(1 - 2k)} \sum_{n=1}^{\infty} n^{2k-1} \frac{q^n}{1 - q^n}.
\]  
Denote by  $q_n := q^n$ and $G(x) := x/(1 - x)$. Applying (2.1), we may rewrite the above expansion as  
\begin{equation}\label{Eisen}
\mathcal{E}_{2k}(z) = 1 + \frac{2}{\zeta(1 - 2k)} \sum_{n=1}^{\infty} n^{2k-1} G(q_n).
\end{equation}

Now we come to the situation of function fields. Let $\Omega=\mathbb{C}_{\infty}\setminus K_{\infty}$ be Drinfeld's  upper half-plane,
which is an analogue of the classical upper half-plane $\mathbb{H}=\{\tau\in\mathbb{C}~ |~ \textrm{Im}(\tau)>0\}.$
Let \( \Gamma = \rm{GL}_2(A) \), \( k \) be a positive integer and \( m \) be an integer with \( 0 \leq m < q - 1 \). In the function fields setting, in analogy with the classical \( \rm{SL}_2(\mathbb{Z}) \)-theory, Goss and Gekeler introduced the modular forms of weight \( k \) and type \( m \) on $\Omega$ associated with \( \Gamma \), named the Drinfeld modular forms (see \cite{Go1,Go2,Gekeler}). The space of such forms is denoted by \( M_{k,m} = M_{k,m}(\Gamma) \), and the subspaces of cusp forms is denoted by  $S_{k,m} = S_{k,m}(\Gamma)$.

Let
\begin{equation}\label{u}
u := u(z) = \frac{1}{e_C(\tilde{\pi} z)}. 
\end{equation}
It is regular on \( \Omega \) and   is the usual uniformizer at $\infty$. For \( a \in A_+ \) of degree \( d \), let
\begin{equation}\label{ua}
u_a := u(az)=\frac{1}{e_C(\tilde{\pi}az)}.  
\end{equation}

In function fields, similarly, we define the Eisenstein series of weight \(k\) with \(k \equiv 0 \pmod{q - 1}\) by
$$
\mathcal{G}_k(z) := \sum \frac{1}{(az + b)^{k}},
$$
where the sum is taken over all nonzero pairs \((a, b) \in A \times A\).
As shown in \cite{Gekeler}, this function admits a decomposition involving its constant term and an \(A\)-expansion. Specifically, we have
\begin{equation}\label{Gekeler}
-\frac{\mathcal{G}_k(z)}{\tilde{\pi}^k} = \sum_{a \in A_+} (\tilde{\pi} a)^{-k} + \sum_{a \in A_+} G_{k,\Lambda_{C}}(u_a),
\end{equation}
where $G_{k, \Lambda_{C}}(z)$ denotes the Goss polynomial associate with the Carlitz module $C$ (see \cite[p. 332, Definition 1]{Goss3}).

This identity illustrates that the Eisenstein series $\mathcal{G}_k(z)$ can be expressed as a sum of a constant term (a ``zeta-value'' at a positive integer) and a series over monic polynomials \(a \in A_+\), each term involving the function $G_{k,\Lambda_{C}}(u_a)$. This structure serves as a prototype for the more general \(A\)-expansions introduced by Petrov.

In \cite{Petrov}, Petrov built on the above example and constructed families of cusp forms with \(A\)-expansions.

\begin{theorem}[{Petrov, see \cite[Theorem 1.3]{Petrov} or \cite[Theorem 1]{Goss3}}]
Let \(k, n\) be two positive integers such that \(k - 2n\) is a positive multiple of \(q - 1\) and \(n \leq p^{\operatorname{ord}_p(k-n)}\). Then
\begin{equation}\label{Petrov}
f_{k,n} := \sum_{a \in A_+} a^{k-n} G_{n, \Lambda_{C}}(u_a)
\end{equation}
is an element of \(S_{k,m}\) with \(n \equiv m \pmod{q-1}\).
\end{theorem}

The above construction led Goss to introduce $\mathfrak{v}$-adic modular forms in function fields in the sense of Serre (see \cite[p.~334, Definition~4]{Goss3}), while Serre \cite{Serre} started to define the $p$-adic modular forms in number fields
as $p$-adic $q$-expansion which are uniform limits of $q$-expansions of the classical Eisenstein series. 

Although similar expansions of Eisenstein series exist in the function field setting due to Gekeler and Petrov (see equations (\ref{Gekeler}) and (\ref{Petrov})), the situation is quite different from that in number fields. This is because the corresponding term $G_{n, \Lambda_{C}}(u_a)$ does not possess a transformation formula analogous to that of the classical Lambert series (\ref{Apostol}) (see \cite{Hamahata}). For the same reason, the transformation formula of Eisenstein series and its expansions (\ref{Gekeler}) may not be used to derive the main result (Theorem \ref{main}).

In fact, the transformation formulae for the logarithm of the Dedekind eta function in number fields and for the Dedekind sum in function fields were discussed in \cite{Hamahata1, Hamahata}. These works also define a new Lambert series in function fields as an analogue of these transformation laws (see (\ref{Def 2.1B})). 

Furthermore, as remarked by Goss \cite{Goss3}:
    ``in a certain sense, Petrov's construction of cusp forms via $A$-expansions is both analogous to and `orthogonal to' (i.e., it lies in the space of cusp forms, rather than Eisenstein series) the classical construction of Eisenstein series'' (see \cite[p.~334, Remark~3]{Goss3}). In the sense of this analogy, comparing (\ref{Eisen}) and   (\ref{Petrov}), we observe that the term $\zeta(1-2k)$ appears in the Lambert  expansions of the classical Eisenstein series, whereas the $A$-expansion (\ref{Petrov}) does not contain any term analogous to $\zeta(1-2k)$. Therefore, we believe it is reasonable that although the term $\zeta(2m+1)$ appears in the classical Ramanujan's identity (\ref{(1.2)}), it does not appear in its function field analogue (\ref{1.7}).

\subsection{Hamahata's Lambert series}
In the classical setting, for any $\tau \in \mathbb{H} = \{ \tau \in \mathbb{C} \mid \operatorname{Im}(\tau) > 0 \}$, the Dedekind eta function is defined as
\[
\eta(\tau) = e^{\pi i \tau / 12} \prod_{n=1}^{\infty} \left(1 - e^{2\pi i n \tau} \right).
\]
Dedekind proved the transformation formula for $\log \eta(\tau)$, that is, for a matrix $\begin{pmatrix} a & b \\ c & d \end{pmatrix} \in \mathrm{SL}_2(\mathbb{Z})$ with $a \neq 0$ and $c > 0$, we have
\begin{equation}\label{eta}
\log \eta\left( \frac{a z + b}{c z + d} \right) = \log \eta(z) + \frac{1}{2} \log\left( \frac{c z + d}{i} \right) + \frac{\pi i (a + d)}{12c} - \pi i D(a, c),
\end{equation}
where $D(a, c)$ is the classical Dedekind sum, defined for coprime integers $a$ and $c > 0$ by

\[
D(a, c) = \frac{1}{4c} \sum_{k=1}^{c-1} \cot\left( \frac{\pi a k}{c} \right) \cot\left( \frac{\pi k}{c} \right).
\]
Since by (\ref{Gn}) we have $$\log \eta(\tau) = \frac{\pi i \tau}{12} - G_1(e^{2\pi i \tau}),$$ the transformation formula for the Lambert series $G_1(e^{2\pi i \tau})$ can be derived from (\ref{eta}).

For the general Lambert series $G_n(e^{2\pi i \tau})$ ($n > 1$), in 1950 Apostol \cite{Apostol} proved   their  transformation formulas, that is, for odd integer $n > 1$  and $\gamma = \begin{pmatrix} a & b \\ c & d \end{pmatrix} \in \mathrm{SL}_2(\mathbb{Z})$, we have
\begin{equation}\label{Apostol}
\begin{aligned}
(c\tau + d)^{n-1} G_n\left( e^{2\pi i \gamma \tau} \right)&= G_n\left( e^{2\pi i \tau} \right) + \frac{1}{2} \zeta(n) \left\{ 1 - (c\tau + d)^{n-1} \right\}\\
&\quad - \frac{(2\pi i)^n}{2(n+1)!} \sum_{r=0}^{n+1} \binom{n+1}{r} (c\tau + d)^{n-r} D_r(-a, c),
\end{aligned}
\end{equation}
where 
\[
D_r(a, c) = \sum_{\mu=1}^{c-1} \overline{B}_{n+1-r}\left( \frac{a \mu}{c} \right) \overline{B}_r\left( \frac{\mu}{c} \right)
\]
is the generalized Dedekind sum and $\overline{B}_n(x)$ be the periodic Bernoulli functions.

In the function fields situation, as remarked in the previous subsection, although the $A$-expansion (\ref{Gekeler}) and (\ref{Petrov}) can be viewed as an analogy of the Fourier expansion for Eisenstein series (\ref{Eisen}), it is not reasonable to view the corresponding terms $G_{n, \Lambda_{C}}(u_a)$ as an analogy of Lambert series, because it does not have a transformation formula similar to (\ref{Apostol}). To fill this gap, Hamahata \cite[Definition~2.1]{Hamahata} introduced the following new series $\xi_{n}(z)$ in the function fields setting.
\begin{definition}[{Hamahata, see \cite[p. 275, Definition 2.1]{Hamahata}}]\label{Def 2.1B}
For a natural number \( n \in \mathbb{N} \), let
\begin{equation}\label{Def 2.1B}
\xi_n(z) = \sum_{a \in A_+} \frac{1}{a^n} u(az).
\end{equation}
\end{definition}
And  he also proved an analogy of Apostol's transformation formula (\ref{Apostol}):
\begin{theorem}[{Hamahata, see \cite[p. 281, Theorem 4.4]{Hamahata}}]
Assume that $q - 1 \mid n + 1$. For $\gamma = \begin{pmatrix} a & b \\ c & d \end{pmatrix} \in GL_2(A)$ with $a \neq 0$, $c \neq 0$, we have
\begin{multline}
(cz + d)^{n-1} \xi_n(\gamma z) = (\det \gamma)^{-1} \Bigg[ \xi_n(z) + \frac{\alpha(n+1)}{\pi c^n} \left\{ (cz + d)^n + \frac{1}{cz + d} \right\} \\
- \frac{1}{\pi c^n} \sum_{r=1}^{n} \alpha(r) \alpha(n+1-r) (cz + d)^{r-1} \\
- \frac{\pi^n}{c^n} \sum_{r=1}^{n} (\det \gamma)^{r-1} (cz + d)^{n-r} t_r^{(n)}(a, c) \Bigg].
\end{multline}
\end{theorem}
Here $t_r^{(n)}(a, c)$ denotes the  function fields analogue of Carlitz's generalized Dedekind sums:
\begin{equation}
t_{r}^{(n)}(a,c)=\sum_{0\neq\mu\in A/cA}e_{C}(\tilde\pi\mu/c)^{-n-1+r}e_{C}(\tilde\pi a\mu/c)^{-r}
\end{equation}
for $a,c\in A\setminus\{0\}$ and $r\in\{0,1,\ldots,n+1\}$.
Combining  (\ref{u}), (\ref{ua}) and (\ref{Def 2.1B}), by noticing that $ \alpha\beta = \tilde{\pi}^2 $, we can state Theorem \ref{main} in the form of Hamahata's Lambert series $\xi_{n}(z)$, which is an analogy of the classical formula (\ref{(1.2+)}) above.
\begin{corollary} Let $\xi_{n}(z)$ be Hamahata's Lambert series (\ref{Def 2.1B}), $\Gamma_{m}$ be Carlitz's  factorial (\ref{factorial}) and $BC_{m}$ be the Bernoulli-Carlitz numbers (\ref{(1.5)}). If $ \alpha $ and $ \beta $ are elements in $\mathbb{C}_{\infty}$ such that $ \alpha\beta = \tilde{\pi}^2 $ and if $ j $ is a positive integer, then we have
\begin{align}
&\quad\alpha^{-(\frac{q-1}{2})(j-1)-(q-2)}\xi_{(q-1)j+(q-2)}\left(\frac{\alpha}{\tilde{\pi}}\right)+\beta^{-(\frac{q-1}{2})(j-1)-(q-2)}\xi_{(q-1)j+(q-2)}\left(\frac{\beta}{\tilde{\pi}}\right)\nonumber\\
&=\sum_{k=0}^{j+1}\frac{BC_{(q-1)k}}{\Gamma_{(q-1)k}}\frac{BC_{(q-1)(j-k+1)}}{\Gamma_{(q-1)(j-k+1)}}\alpha^{(\frac{q-1}{2})k}\beta^{(\frac{q-1}{2})(j-k+1)}. \label{1.7+}
\end{align}
\end{corollary}
As demonstrated in \cite{DKK}, our result also opens the possibility of placing Ramanujan's identity within a broader analytic framework in the function field setting.

\section{Proof of Theorem \ref{main}}
Our proof is inspired by Chavan \cite{Chavan} for his elementary proof of Ramanujan's formula (\ref{(1.2)}).

By (\ref{(1.6)}), we get
\begin{equation}
\begin{aligned}
e_{C}(z)&=e_{\tilde{\pi}A}(z) \\
&=z\prod_{0\neq \lambda\in\tilde{\pi}A}\left(1-\frac{z}{\lambda}\right)\\
&=z\prod_{0\neq a\in A}\left(1-\frac{z}{\tilde{\pi}a}\right)\\
&=z\prod_{a\in A_{+}}\left(1-\left(\frac{z}{\tilde{\pi}a}\right)^{q-1}\right),
\end{aligned}
\end{equation}
the last identity follows from the polynomial factorization
$$x^{q-1}-y^{q-1}=\prod_{\xi\in\mathbb{F}_{q}^{*}}(x-\xi y).$$
From (\ref{1.6+}),  
$$e_{C}'(z)=1.$$
So
\begin{equation}
\begin{aligned}
\frac{1}{e_{C}(z)}&=\frac{e_{C}'(z)}{e_{C}(z)} \\
&=\frac{1}{z}+\sum_{a\in A_{+}}\frac{-\frac{(q-1)z^{q-2}}{(\tilde{\pi}a)^{q-1}}}{1-\left(\frac{z}{\tilde{\pi}a}\right)^{q-1}}\\
&=\frac{1}{z}-z^{q-2}\sum_{a\in A_{+}}\frac{1}{z^{q-1}-(\tilde{\pi}a)^{q-1}}.
\end{aligned}
\end{equation}
Thus for any $\alpha\in\Omega$ and  $b\in A$, we have
\begin{equation}\label{(1.7)}
\frac{1}{e_{C}(\alpha b)}=\frac{1}{\alpha b}-(\alpha b)^{q-2}\sum_{a\in A_{+}}\frac{1}{(\alpha b)^{q-1}-(\tilde{\pi}a)^{q-1}}.
\end{equation}
Note that, by (\ref{(1.7)}) and (\ref{(1.3)}), for any positive integer $j,$ we have
\begin{align*}
&\quad\alpha^{(\frac{q-1}{2})(j+1)}\sum_{b\in A_{+}}\sum_{a\in A_{+}}\frac{1}{(\alpha b)^{(q-1)j}\left((\alpha b)^{q-1}-(\tilde{\pi}a)^{q-1}\right)}\\
&=\alpha^{(\frac{q-1}{2})(j+1)}\sum_{b\in A_{+}}\frac{1}{(\alpha b)^{(q-1)j}}\frac{1}{(\alpha b)^{q-2}}\left(\frac{1}{\alpha b}-\frac{1}{e_{C}(\alpha b)}\right)\\
&=\alpha^{(\frac{q-1}{2})(j+1)}\left(\frac{1}{\alpha^{(q-1)(j+1)}}\zeta_{\infty}\left((q-1)(j+1)\right) \right. \\
&\left.\qquad\qquad\qquad\qquad\qquad\quad -\sum_{b\in A_{+}}\frac{1}{(\alpha b)^{(q-1)j+(q-2)}e_{C}(\alpha b)}\right).
\end{align*}
Then by applying (\ref{(1.4)}) and notice that $\alpha\beta=\tilde{\pi}^{2}$, we further get
\begin{equation}\label{(2.4+)}
\begin{aligned}
&\quad\alpha^{(\frac{q-1}{2})(j+1)}\sum_{b\in A_{+}}\sum_{a\in A_{+}}\frac{1}{(\alpha b)^{(q-1)j}\left((\alpha b)^{q-1}-(\tilde{\pi}a)^{q-1}\right)}\\
&=\alpha^{(\frac{q-1}{2})(j+1)}\frac{1}{\alpha^{(q-1)(j+1)}}\tilde{\pi}^{(q-1)(j+1)} \frac{BC_{(q-1)(j+1)}}{\Gamma_{(q-1)(j+1)}}\\
&\quad-\alpha^{(\frac{q-1}{2})(j+1)}\sum_{b\in A_{+}}\frac{1}{(\alpha b)^{(q-1)j+(q-2)}e_{C}(\alpha b)}\\
&=\beta^{(\frac{q-1}{2})(j+1)}\frac{BC_{(q-1)(j+1)}}{\Gamma_{(q-1)(j+1)}}-\alpha^{-(\frac{q-1}{2})(j-1)-(q-2)} \\
&\quad\times\sum_{b\in A_{+}}\frac{b^{-(q-1)j-(q-2)}}{e_{C}(\alpha b)}.
\end{aligned}
\end{equation}
On the other hand, a direct calculation shows that
\begin{align*}
&\quad\sum_{b\in A_{+}}\sum_{a\in A_{+}}\frac{1}{(\alpha b)^{(q-1)j}\left((\alpha b)^{q-1}-(\tilde{\pi}a)^{q-1}\right)}\\
&=-\sum_{b\in A_{+}}\frac{1}{(\alpha b)^{(q-1)j}}\sum_{a\in A_{+}}\frac{1}{(\tilde{\pi}a)^{q-1}}\\
&\quad+\sum_{b\in A_{+}}\sum_{a\in A_{+}}\frac{1}{(\alpha b)^{(q-1)(j-1)}(\tilde{\pi}a)^{q-1}\left((\alpha b)^{q-1}-(\tilde{\pi}a)^{q-1}\right)}\\
&=-\sum_{b\in A_{+}}\frac{1}{(\alpha b)^{(q-1)j}}\sum_{a\in A_{+}}\frac{1}{(\tilde{\pi}a)^{q-1}} \\
&\quad-\sum_{b\in A_{+}}\frac{1}{(\alpha b)^{(q-1)(j-1)}}\sum_{a\in A_{+}}\frac{1}{(\tilde{\pi}a)^{(q-1)2}}\\
&\quad+\sum_{b\in A_{+}}\sum_{a\in A_{+}}\frac{1}{(\alpha b)^{(q-1)(j-2)}(\tilde{\pi}a)^{(q-1)2}\left((\alpha b)^{q-1}-(\tilde{\pi}a)^{q-1}\right)}\\
&=\cdots\\
&\quad (\textrm{inductively})\\
&=-\sum_{b\in A_{+}}\frac{1}{(\alpha b)^{(q-1)j}}\sum_{a\in A_{+}}\frac{1}{(\tilde{\pi}a)^{q-1}} \\
&\quad-\sum_{b\in A_{+}}\frac{1}{(\alpha b)^{(q-1)(j-1)}}\sum_{a\in A_{+}}\frac{1}{(\tilde{\pi}a)^{(q-1)2}}\\
&\quad\cdots\\
&\quad-\sum_{b\in A_{+}}\frac{1}{(\alpha b)^{q-1}}\sum_{a\in A_{+}}\frac{1}{(\tilde{\pi}a)^{(q-1)j}} \\
&\quad+\sum_{b\in A_{+}}\sum_{a\in A_{+}}\frac{1}{(\tilde{\pi}a)^{(q-1)j}}\frac{1}{(\alpha b)^{q-1}-(\tilde{\pi}a)^{q-1}}\\
&=-\sum_{k=1}^{j}\left(\sum_{b\in A_{+}}\frac{1}{(\alpha b)^{(q-1)(j-k+1)}}\sum_{a\in A_{+}}\frac{1}{(\tilde{\pi}a)^{(q-1)k}}\right)\\
&\quad+\sum_{b\in A_{+}}\sum_{a\in A_{+}}\frac{1}{(\tilde{\pi}a)^{(q-1)j}}\frac{1}{(\alpha b)^{q-1}-(\tilde{\pi}a)^{q-1}}\\
&=-\sum_{k=1}^{j}\frac{\zeta_{\infty}\left((q-1)(j-k+1)\right)}{\alpha^{(q-1)(j-k+1)}}\frac{\zeta_{\infty}\left((q-1)k\right)}{\tilde{\pi}^{(q-1)k}}\\
&\quad+\sum_{b\in A_{+}}\sum_{a\in A_{+}}\frac{1}{(\tilde{\pi}a)^{(q-1)j}}\frac{1}{(\alpha b)^{q-1}-(\tilde{\pi}a)^{q-1}},
\end{align*}
the last identity follows from (\ref{(1.3)}).

Then multiplying both sides of the above equation by $\alpha^{(\frac{q-1}{2})(j+1)},$ notice (\ref{(1.4)}) and the relation
$\alpha\beta=\tilde{\pi}^{2}$, we have
\begin{equation}
\begin{aligned}
&\quad\alpha^{(\frac{q-1}{2})(j+1)}\sum_{b\in A_{+}}\sum_{a\in A_{+}}\frac{1}{(\alpha b)^{(q-1)(j+1)}\left((\alpha b)^{q-1}-(\tilde{\pi}a)^{q-1}\right)}\\
&=-\sum_{k=1}^{j}\frac{\zeta_{\infty}\left((q-1)(j-k+1)\right)}{\tilde{\pi}^{(q-1)(j-k+1)}}\frac{\zeta_{\infty}\left((q-1)k\right)}{\tilde{\pi}^{(q-1)k}}\alpha^{(\frac{q-1}{2})k}\beta^{(\frac{q-1}{2})(j-k+1)}\\
&\quad+\alpha^{(\frac{q-1}{2})(j+1)}\sum_{b\in A_{+}}\sum_{a\in A_{+}}\frac{1}{(\tilde{\pi}a)^{(q-1)j}}\frac{1}{(\alpha b)^{q-1}-(\tilde{\pi}a)^{q-1}}\\
&=-\sum_{k=1}^{j}\frac{BC_{(q-1)(j-k+1)}}{\Gamma_{(q-1)(j-k+1)}}\frac{BC_{(q-1)k}}{\Gamma_{(q-1)k}}\alpha^{(\frac{q-1}{2})k}\beta^{(\frac{q-1}{2})(j-k+1)}\\
&\quad+\alpha^{(\frac{q-1}{2})(j+1)}\sum_{b\in A_{+}}\sum_{a\in A_{+}}\frac{1}{(\tilde{\pi}a)^{(q-1)j}}\frac{1}{(\alpha b)^{q-1}-(\tilde{\pi}a)^{q-1}}.
\end{aligned}
\end{equation}
By comparing with (\ref{(2.4+)}) we get
\begin{equation}\label{(1.9)}
\begin{aligned}
&\quad-\alpha^{-(\frac{q-1}{2})(j-1)-(q-2)}\sum_{b\in A_{+}}\frac{b^{-(q-1)j-(q-2)}}{e_{C}(\alpha b)}\\
&=-\sum_{k=0}^{j}\frac{BC_{(q-1)(j-k+1)}}{\Gamma_{(q-1)(j-k+1)}}\frac{BC_{(q-1)k}}{\Gamma_{(q-1)k}}\alpha^{(\frac{q-1}{2})k}\beta^{(\frac{q-1}{2})(j-k+1)}\\
&\quad+\alpha^{(\frac{q-1}{2})(j+1)}\sum_{b\in A_{+}}\sum_{a\in A_{+}}\frac{1}{(\tilde{\pi}a)^{(q-1)j}}\frac{1}{(\alpha b)^{q-1}-(\tilde{\pi}a)^{q-1}}.
\end{aligned}
\end{equation}
To compute the last term, note that $\alpha\beta=\tilde{\pi}^{2}$, we have
\begin{align*}
&\quad\sum_{b\in A_{+}}\sum_{a\in A_{+}}\frac{1}{(\tilde{\pi}a)^{(q-1)j}}\frac{1}{(\alpha b)^{q-1}-(\tilde{\pi}a)^{q-1}}\\
&=\frac{1}{\tilde{\pi}^{(q-1)j}}\sum_{b\in A_{+}}\sum_{a\in A_{+}}\frac{1}{a^{(q-1)j}}\frac{1}{(\alpha b)^{q-1}-\alpha^{\frac{q-1}{2}}\beta^{\frac{q-1}{2}}a^{q-1}}\\
&=\frac{1}{\tilde{\pi}^{(q-1)j}}\sum_{b\in A_{+}}\sum_{a\in A_{+}}\frac{1}{\alpha^{\frac{q-1}{2}}a^{(q-1)j}}\frac{1}{\alpha^{\frac{q-1}{2}}b^{q-1}-\beta^{\frac{q-1}{2}}a^{q-1}}\\
&=\frac{1}{\alpha^{(\frac{q-1}{2})j}\beta^{(\frac{q-1}{2})j}}\sum_{a\in A_{+}}\frac{\beta^{\frac{q-1}{2}}}{\alpha^{\frac{q-1}{2}}}\frac{1}{a^{(q-1)j}}\frac{1}{\tilde{\pi}^{q-1}b^{q-1}-\beta^{q-1}a^{q-1}}\\
&=-\frac{\beta^{(\frac{q-1}{2})j}}{\alpha^{(\frac{q-1}{2})j}}\frac{\beta^{\frac{q-1}{2}}}{\alpha^{\frac{q-1}{2}}}\sum_{a\in A_{+}}\frac{1}{(\beta a)^{(q-1)j}}\frac{1}{(\beta a)^{q-1}-({\tilde{\pi}b})^{q-1}}.
\end{align*}
Then multiplying both sides of the above equation by $\alpha^{(\frac{q-1}{2})(j+1)}$, notice (\ref{(1.7)}) and (\ref{(1.3)}), we get
\begin{equation}\label{(1.10)}
\begin{aligned}
&\quad\alpha^{(\frac{q-1}{2})(j+1)}\sum_{b\in A_{+}}\sum_{a\in A_{+}}\frac{1}{(\tilde{\pi}a)^{(q-1)j}}\frac{1}{(\alpha b)^{q-1}-(\tilde{\pi}a)^{q-1}}\\
&=-\beta^{(\frac{q-1}{2})(j+1)}\sum_{a\in A_{+}}\frac{1}{(\beta a)^{(q-1)j}}\frac{1}{(\beta a)^{q-2}}\left(\frac{1}{\beta a}-\frac{1}{e_{C}(\beta a)}\right)\\
&=-\beta^{(\frac{q-1}{2})(j+1)}\left(\frac{1}{\beta^{(q-1)(j+1)}}\zeta_{\infty}\left((q-1)(j+1)\right) \right. \\
&\qquad\qquad\qquad\qquad\qquad\quad -\left.\sum_{a\in A_{+}}\frac{1}{(\beta a)^{(q-1)j+(q-2)}e_{C}(\beta a)}\right).
\end{aligned}
\end{equation}
Similarly with the reasoning of (\ref{(2.4+)}), by applying (\ref{(1.4)}) to the above equation and notice that $\alpha\beta=\tilde{\pi}^{2}$, we see that
\begin{equation}
\begin{aligned}
&\quad\alpha^{(\frac{q-1}{2})(j+1)}\sum_{b\in A_{+}}\sum_{a\in A_{+}}\frac{1}{(\tilde{\pi}a)^{(q-1)j}}\frac{1}{(\alpha b)^{q-1}-(\tilde{\pi}a)^{q-1}}\\
&=-\beta^{(\frac{q-1}{2})(j+1)}\frac{1}{\beta^{(q-1)(j+1)}}\tilde{\pi}^{(q-1)(j+1)} \frac{BC_{(q-1)(j+1)}}{\Gamma_{(q-1)(j+1)}}\\
&\quad+\beta^{(\frac{q-1}{2})(j+1)}\sum_{a\in A_{+}}\frac{1}{(\beta a)^{(q-1)j+(q-2)}e_{C}(\beta a)}\\
&=-\alpha^{(\frac{q-1}{2})(j+1)}\frac{BC_{(q-1)(j+1)}}{\Gamma_{(q-1)(j+1)}} \\
&\quad+\beta^{-(\frac{q-1}{2})(j-1)-(q-2)}\sum_{a\in A_{+}}\frac{a^{-(q-1)j-(q-2)}}{e_{C}(\beta a)}.
\end{aligned}
\end{equation}
Finally, by substituting the above equation into (\ref{(1.9)}) we obtain
\begin{align*}
&\alpha^{-(\frac{q-1}{2})(j-1)-(q-2)}\sum_{b\in A_{+}}\frac{b^{-(q-1)j-(q-2)}}{e_{C}(\alpha b)} \\
&\quad +\beta^{-(\frac{q-1}{2})(j-1)-(q-2)}\sum_{a\in A_{+}}\frac{a^{-(q-1)j-(q-2)}}{e_{C}(\beta a)}\\
&=\sum_{k=0}^{j+1}\frac{BC_{(q-1)(j-k+1)}}{\Gamma_{(q-1)(j-k+1)}}\frac{BC_{(q-1)k}}{\Gamma_{(q-1)k}}\alpha^{(\frac{q-1}{2})k}\beta^{(\frac{q-1}{2})(j-k+1)},
\end{align*}
which is the desired result.

\section*{Acknowledgements} 
 The authors are enormously grateful to the anonymous referee for his/her careful
reading of this paper, and for his/her insightful comments and suggestions. 

Su Hu is supported by the Natural Science Foundation of Guangdong Province, China (No. 2024A1515012337).  Min-Soo Kim is supported by the National Research Foundation of Korea(NRF) grant funded by the Korea government(MSIT) (No. NRF-2022R1F1A1065551 and No. RS-2025-23323090). 

\section*{Statements and Declarations}
\subsection*{Declaration of competing interest}
The authors declare that they have no known competing financial interests or personal relationships that could have appeared to influence the work reported in this paper.

Authors declare that they do not have any conflict of interest.

\subsection*{Data availability}
No data was used for the research described in the article.

\end{document}